\documentclass{article}

\usepackage{amsmath,amssymb,amsfonts}
\usepackage[T2A]{fontenc}
\usepackage[cp1251]{inputenc}
\usepackage[russian,english]{babel}
\usepackage[12pt]{extsizes}
\newcommand{\7}{\!\!\!\!\!\!\!}
\newcommand{\5}{\!\!\!\!\!}

\begin{document}
\begin{center}
\textbf{Oscillation of irrational measure function in the multidimensional case.\\
D. O. Shatskov}
\end{center}
\begin{abstract}
Let $\Theta$ is a matrix of size $m\times n$. We denote the irrational measure function, as
$$
\psi_\Theta (t)=\min_{\substack{
  x_i{\in}\mathbb{Z}\\
  1\leqslant\max\limits_{1\leqslant i\leqslant n}|x_i|\leqslant t}} \max\limits_{1\leqslant j\leqslant
m}\|\theta^1_jx_1+\ldots+\theta^n_jx_n\|.
$$
We proved that difference function $\psi_\Theta-\psi_{\Theta'}$ for almost all pairs $\Theta$, $\Theta'$ in cases $m=1$, $n=2$ or
$m\geqslant2$ and $n=1$ changes its sign infinity many times as $t\rightarrow+\infty$.

Bibliography: 5 titles.
\end{abstract}

\footnotetext{Research is supported by the grant RNF 14-11-00433.}%

For a real matrix $$\Theta=\left(%
\begin{array}{cccc}
  \theta^1_1 &  \ldots & \theta^n_1\\
  \vdots &  \ddots&  \vdots\\
  \theta^1_m &  \ldots & \theta^n_m\\
\end{array}%
\right)$$ we consider the irrationality measure function
$$
\psi_\Theta (t)=\min_{\substack{
  x_i{\in}\mathbb{Z}\\
  1\leqslant\max\limits_{1\leqslant i\leqslant n}|x_i|\leqslant t\\
  }}
\max\limits_{1\leqslant j\leqslant
m}\|\theta^1_jx_1+\ldots+\theta^n_jx_n\|,
$$
here $\|\cdot\|$ stands for the distance to the nearest integer. For two matrix $\Theta$ and $\Theta'$ we consider the  difference function $\psi_\Theta
(t) - \psi_{\Theta'}(t)$.

If $m=n=1$ then $\Theta=\alpha$ and  
$$
\psi_\alpha(t)=\min\limits_{1\leqslant q\leqslant t}\|q\alpha\|.
$$

In the paper  \cite{MosKan} it is proved that for any two real numbers $\alpha$ and $\beta$, such that $\alpha \pm \beta {\not\in} \mathbb{Z}$ the difference function
$$
\psi_\alpha (t) - \psi_\beta (t)
$$
changes its sign infinitely many often as $t\rightarrow+\infty$.

Two theorems below are well-known.
 
{\bf Theorem 1.} (A. Khintchine, 1926). Let a function $\psi(t)$ decrease to zero as $t\rightarrow+\infty$. Then there exist two algebraically independent real numbers $\alpha$ and $\beta$ such that for all $t$ large enough one has
$$
\psi_{(\alpha\,\,\beta)}(t)=\min\limits_{1\leqslant\max(|x_1|,|x_2|)\leqslant
t}\|x_1\alpha+x_2\beta\|\leqslant\psi(t).
$$

{\bf Theorem 2.} (A. Khintchine, 1926). Let a function $\psi(t)$ decrease to zero as $t\rightarrow+\infty$ and the function $t\psi(t)$ increase to infinity as $t\rightarrow+\infty$. Then there exist two algebraically independent real numbers $\alpha$ and $\beta$
such that for all $t$ large enough
$$
\psi_{\scriptsize\left(
\begin{array}{c}
  \alpha\\
  \beta\\
\end{array}
\right)}(t)=\min\limits_{1\leqslant x\leqslant
t}\max(\|x\alpha\|,\|x\beta\|)\leqslant\psi(t).
$$

This shows that there is no general oscillation property for difference function in any dimension greater than one.

The main result of this note is the following

{\bf Theorem.} Let $m=1$ and $n=2$ or $m\geqslant2$ and $n=1$,
than for almost all matrix $\Theta$ and $\Theta'$ of size $m\times
n$ the difference function
$$
\psi_\Theta (t) - \psi_{\Theta'}(t)
$$
changes its sign infinitely many times as $t\rightarrow+\infty$.

\newpage

\begin{center}
\textbf{Осцилляция функции меры иррациональности в многомерном случае\\
Д.О. Шацков}
\end{center}

\footnotetext{Работа выполнена при поддержке гранта РНФ 14-11-00433.}%

\textbf{1. Введение.}
Для матрицы $$\Theta=\left(%
\begin{array}{cccc}
  \theta^1_1 &  \ldots & \theta^n_1\\
  \vdots &  \ddots&  \vdots\\
  \theta^1_m &  \ldots & \theta^n_m\\
\end{array}%
\right)$$ рассмотрим функцию меры иррациональности

$$
\psi_\Theta (t)=\min_{\substack{
  x_i{\in}\mathbb{Z}\\
  1\leqslant\max\limits_{1\leqslant i\leqslant n}|x_i|\leqslant t}}
\max\limits_{1\leqslant j\leqslant
m}\|\theta^1_jx_1+\ldots+\theta^n_jx_n\|,
$$
где $\|\cdot\|$ обозначает расстояние до ближайшего целого. Для
двух матриц $\Theta$ и $\Theta'$ составим разность $\psi_\Theta
(t) - \psi_{\Theta'}(t)$. Нас интересует как ведет себя эта
разность при $t\rightarrow+\infty$.

В одномерном случае $\Theta=\alpha$ и функцию $\psi_\Theta$ можно
записать так
$$
\psi_\alpha(t)=\min\limits_{1\leqslant q\leqslant t}\|q\alpha\|.
$$

Н.Г. Мощевитин и И.Д. Кан в совместной работе  \cite{MosKan}
доказали, что для двух вещественных чисел $\alpha$ и $\beta$,
таких что $\alpha \pm \beta {\not\in} \mathbb{Z}$ разность
$$
\psi_\alpha (t) - \psi_\beta (t)
$$
бесконечно много раз меняет знак при $t\rightarrow+\infty$.

Известны две теоремы.

{\bf Теорема} (Хинчин, 1926). Пусть функция $\psi(t)$ убывает к
нулю при $t\rightarrow+\infty$. Тогда существуют два алгебраически
независимых числа $\alpha$ и $\beta$ такие, что для достаточно
больших $t$ выполняется
$$
\psi_{(\alpha\,\,\beta)}(t)=\min\limits_{1\leqslant\max(|x_1|,|x_2|)\leqslant
t}\|x_1\alpha+x_2\beta\|\leqslant\psi(t).
$$

{\bf Теорема} (Хинчин, 1926). Пусть функция $\psi(t)$ убывает к
нулю при $t\rightarrow+\infty$ и функция $t\psi(t)$ бесконечно
возрастает при $t\rightarrow+\infty$. Тогда существуют два
алгебраически независимых действительных числа $\alpha$ и $\beta$
таких, что для достаточно больших $t$ выполняется
$$
\psi_{\scriptsize\left(
\begin{array}{c}
  \alpha\\
  \beta\\
\end{array}
\right)}(t)=\min\limits_{1\leqslant x\leqslant
t}\max(\|x\alpha\|,\|x\beta\|)\leqslant\psi(t).
$$

Следствием этих двух теорем служит утверждение, что для больших
размерностей в общем случае осцилляции нет.

В данной работе будет доказан следующий результат

{\bf Теорема.} Пусть $m=1$ и $n=2$ или $m\geqslant2$ и $n=1$,
тогда для почти всех матриц $\Theta$ и $\Theta'$ размера $m\times
n$ разность
$$
\psi_\Theta (t) - \psi_{\Theta'}(t)
$$
осциллирует бесконечное число раз при $t\rightarrow+\infty$.

\textbf{2. Изучение функции $\psi_\Theta (t)$ в случае $m=1$ и $n=2$.}

Далее в работе $0<\lambda<\lambda_0$, $(\lambda_0$ эффективно
вычисляется из лемм 3, 4, 5, 13 и 14) $0<\varepsilon<1$
$(\varepsilon$ эффективно вычисляется из лемм 3, 4, 5, 12 и 14),
$0<\delta<1$ $(\delta$ эффективно вычисляется из леммы 14) и
$k>K_0=K_0(\lambda)$, ($K_0$ эффективно вычисляется из лемм 2, 3,
4, 5, 8, 9, 11, 13, 14 и 15).

В случае $m=1$ и $n=2$  функция $\psi_\Theta (t)$ имеет вид
$$
\psi_{(\alpha\,\,\beta)}(t)=\min\limits_{1\leqslant\max(|x_1|,|x_2|)\leqslant
t}\|x_1\alpha+x_2\beta\|.
$$
На квадрате $[0;1]^2$ рассмотрим два множества:
$$
\underline{M_k}=\left\{(\alpha,\beta){\in}[0;1]^2:
\psi_{(\alpha\,\,\beta)}(k)>\frac{\varepsilon}{k^2}\right\}=$$
$$
=\left\{(\alpha,\beta){\in}[0;1]^2: \forall
(x_1,x_2){\in}\mathbb{Z}^2 \,\,\,
1\leqslant\max(|x_1|,|x_2|)\leqslant k
\,\,\,\,\|x_1\alpha+x_2\beta\|>\frac{\varepsilon}{k^2}\right\}=$$
$$
=\left\{(\alpha,\beta){\in}[0;1]^2:\forall
(x_1,x_2){\in}\mathbb{Z}^2\,\,
{1\leqslant\max(|x_1|,|x_2|)\leqslant k} \, \forall
q{\in}\mathbb{Z}
\,\,\,|x_1\alpha+x_2\beta-q|>\frac{\varepsilon}{k^2}\right\};
$$
$$
\overline{M}_k=\left\{(\alpha,\beta){\in}[0;1]^2:
\psi_{(\alpha\,\,\beta)}(k)\leqslant\frac{\varepsilon}{k^2}\right\}=$$
$$
=\left\{(\alpha,\beta){\in}[0;1]^2: \exists
(x_1,x_2){\in}\mathbb{Z}^2 \,\,\,
1\leqslant\max(|x_1|,|x_2|)\leqslant k
\,\,\,\,\|x_1\alpha+x_2\beta\|\leqslant\frac{\varepsilon}{k^2}\right\}=$$
$$
=\left\{(\alpha,\beta){\in}[0;1]^2:\exists
(x_1,x_2){\in}\mathbb{Z}^2\,\,\,
1\leqslant\max(|x_1|,|x_2|)\leqslant k\,\, \exists
q{\in}\mathbb{Z} \,\,\,\,
|x_1\alpha+x_2\beta{-}q|\leqslant\frac{\varepsilon}{k^2}\right\}.
$$

Верно равенство
\begin{equation}\label{raven}
\mu(\underline{M_k})=1-\mu(\overline{M}_k).
\end{equation}

Рассмотрим плоскость $O\alpha\beta$. Для любых целых чисел $x_1$,
$x_2$ и $q$ зададим множество
$$
A(x_1,x_2,q)=\left\{(\alpha,\beta){\in}\mathbb{R}^2:|x_1\alpha+x_2\beta-q|\leqslant\frac{\varepsilon}{k^2}\right\}.
$$
Обозначим
$$
A(x_1,x_2)=\left\{(\alpha,\beta){\in}\mathbb{R}^2: \exists
q{\in}\mathbb{Z}\,
|x_1\alpha+x_2\beta-q|\leqslant\frac{\varepsilon}{k^2}\right\}=\bigcup\limits_{q{\in}\,\mathbb{Z}}
A(x_1,x_2,q).
$$

Выпишем свойства множества $A(x_1,x_2)$:
$$
A(x_1,x_2)+\left(\frac{1}{x_1}\mathbb{Z}\right)\times\left(\frac{1}{x_2}\mathbb{Z}\right)=A(x_1,x_2);
$$
$$
A(x_1,x_2)=A(-x_1,-x_2);
$$
$$
\mu\left(A(x_1,x_2){\cap}[0;1]^2\right)=\frac{2\varepsilon}{k^2},
$$
где $\mu(\cdot)$ - мера Лебега, а знак $\times$ обозначает прямое
произведение множеств.

Для любых двух линейно независимых пар $(x_1,x_2)$ и $(y_1,y_2)$ и
любых $q_1$ и $q_2$ пересечением множеств $A(x_1,x_2,q_1)$ и
$A(y_1,y_2,q_2)$ будет параллелограмм с площадью
\begin{equation}\label{plr}
\mu\left(A(x_1,x_2,q_1){\cap}
A(y_1,y_2,q_2)\right)=\frac{1}{\left|\det\left(
\begin{array}{cc}
x_1 & x_2 \\
y_1 & y_2 \\
\end{array}
\right)\right|}\frac{4\varepsilon^2}{k^4}.
\end{equation}

{\bf Замечание 1.} Диаметр этого параллелограмма не более
$\frac{4\varepsilon}{k}$.

Обозначим через $\Lambda$ - решетку, построенную на векторах
$\left(\frac{y_2}{\Delta},\frac{-x_2}{\Delta}\right)$ и
$\left(\frac{-y_1}{\Delta},\frac{x_1}{\Delta}\right)$.

{\bf Замечание 2.} Центры параллелограммов, получаемых при
переборе всевозможных целочисленных значений $q_1$, $q_2$ лежат в
узлах решетки $\Lambda$.

Верно равенство
\begin{equation}\label{pk}
\overline{M}_k=\7\5\bigcup_{\substack{
  {(x_1,x_2)}{\in}\mathbb{Z}^2\\
  1\leqslant\max(|x_1|, |x_2|)\leqslant k}}
  \7\left(A(x_1,x_2){\cap}[0;1]^2\right).
\end{equation}

Для удобства здесь и далее в работе через $S$ будем обозначать
параллелепипеды, соответствующей размерности, со сторонами
параллельными координатным гиперплоскостям.

{\bf Теорема Ярника.} Пусть $G$ --- выпуклая область на плоскости,
$N$ --- число целых точек в области $G$, $P$ --- площадь области
$G$, $L$ --- периметр области $G$, $L\geqslant1$, тогда выполнены
неравенства
$$
P-L<N<L+P.
$$

{\bf Лемма 1.} Пусть $S$ ---
квадрат со стороной $\lambda$, $\Delta=\left|\det\left(%
\begin{array}{cc}
  x_1 & x_2\\
  y_1 & y_2\\
\end{array}%
\right)\right|$, $\Delta\not=0$, $x_1, x_2, y_1,
y_2{\in}\mathbb{Z}$ и $x_1, x_2, y_1, y_2>\frac{1}{\lambda}$,  $N$
- число точек решетки $\Lambda$ в квадрате $S$, тогда выполняется
неравенство
$$
N<\lambda^2\Delta+2\lambda\sqrt{x_1^2+y_1^2}+2\lambda\sqrt{x_2^2+y_2^2}.
$$

{\bf Доказательство.} Подействуем на решетку $\Lambda$ справа линейным преобразованием $T=\left(%
\begin{array}{cc}
  x_1 & x_2\\
  y_1 & y_2\\
\end{array}%
\right)$, которое переводит решетку $\Lambda$ в ортонормированную
решетку~$\mathbb{Z}^2$. Квадрат $S$ при этом преобразовании
перейдет в параллелограмм с площадью $\lambda^2\Delta$, длины
сторон этого параллелограмма будут равны
$\lambda\sqrt{x_1^2+y_1^2}$ и $\lambda\sqrt{x_2^2+y_2^2}$.
Применим теорему Ярника и получим неравенство.

Лемма 1 доказана.

Обозначим
$E_k=\left\{(x_1,x_2){\in}\mathbb{N}^2:(x_1,x_2)=1,\frac{k}{2}\leqslant
x_1,x_2\leqslant k\right\}$.

{\bf Лемма 2.} Пусть $\Delta=\left|\det\left(%
\begin{array}{cc}
  x_1 & x_2\\
  y_1 & y_2\\
\end{array}%
\right)\right|$, тогда выполняется неравенство
\[
\sum_{\substack{
  (x_1,\,\,x_2){\in} E_k\\
  (y_1,\,\,y_2){\in} E_k\\
  (x_1,x_2)\neq(y_1,y_2)\\
  }}
\frac{1}{\Delta}\leqslant 9k^2\ln k.
\]
{\bf Доказательство.} Рассмотрим плоскость $Oy_1y_2$. Для точки
$(x_1,x_2){\in} E_k$ построим семейство прямых $x_1y_2-x_2y_1=q$,
где $q{\in}\mathbb{Z}$. Уравнение прямой из этого семейства можно
переписать в
таком виде $\det\left(%
\begin{array}{cc}
  x_1 & x_2\\
  y_1 & y_2\\
\end{array}\right)=q$. Все точки решетки $\mathbb{Z}^2$ лежат на этих
прямых. Расстояние между соседними прямыми из этого семейства
равно $\frac{1}{\sqrt{x_1^2+x_2^2}}$. Расстояние между соседними
целыми точками на каждой прямой из этого семейства равно
$\sqrt{x_1^2+x_2^2}\geqslant\frac{k}{\sqrt{2}}$. На каждой прямой
будет не более 2 точек $(y_1,y_2){\in} E_k$. Так как мы
рассматриваем пары $(x_1,x_2){\in}E_k$ и $(y_1,y_2){\in}E_k$, то
$\Delta\leqslant k^2$. Оценим интересующую нас сумму
$$
\sum_{\substack{
  (x_1,x_2){\in} E_k\\
  (y_1,y_2){\in} E_k\\
  (x_1,x_2)\neq(y_1,y_2)}}
  \7\frac{1}{\Delta}=\sum\limits_{(x_1,\,x_2){\in}
E_k}\sum\limits_{n=1}^{k^2}\sum_{\substack{
  \Delta=n\\
  (y_1,y_2){\in}E_k\\
  }}\frac{1}{\Delta}
\leqslant\sum\limits_{(x_1,\,x_2){\in} E_k}
\sum\limits_{n=1}^{k^2}\frac{4}{n}\leqslant$$
$$
\leqslant\sum\limits_{(x_1,\,x_2){\in} E_k} \left(4+8\ln
k\right)\leqslant 9k^2\ln k.
$$
Лемма 2 доказана.

{\bf Лемма 3.} Для любого квадрата $S$ со стороной $\lambda$,
верны неравенства
$$
\lambda^2\left(
\frac{\varepsilon}{3\zeta(2)}-37\frac{\varepsilon^2}{\zeta^2(2)}\right)\leqslant\mu(\overline{M}_k{\cap}
S)\leqslant 5\varepsilon\lambda^2.
$$

{\bf Доказательство.} Через $S(x_1,x_2,k_1,k_2)$ прямоугольник со
сторонами $\frac{k_1}{x_1}$ и $\frac{k_2}{x_2}$, где
$x_1,x_2,k_1,k_2{\in}\mathbb{N}_0$. В случае $x_1=0$ или $x_2=0$
под $S$ понимается бесконечная полоса параллельная оси $O\alpha$
или оси $O\beta$ соответственно.

Для любого квадрата $S$ и любой пары чисел $(x_1,x_2)$, если
положить $k_1=\left[\lambda x_1\right]$ и $k_2=\left[\lambda
x_2\right]$, то можно указать два прямоугольника, центры которых
совпадают с центром квадрата $S$ и выполняется условие
$S(x_1,x_2,k_1,k_2){\subset}S{\subset}S(x_1,x_2,k_1+1,k_2+1)$.

Запишем ограничение снизу на меру пересечения $\overline{M}_k$ и
$S$
\begin{equation}\label{osn}
\mu\left(\overline{M}_k{\cap}S\right)\geqslant\5\!\sum\limits_{(x_1,\,x_2){\in}
E_k}\!\!\!\mu(A(x_1,x_2){\cap}
S)-\7\sum\limits_{\substack{
  (x_1,\,\,x_2){\in} E_k\\
  (y_1,\,\,y_2){\in} E_k\\
  (x_1,x_2)\neq(y_1,y_2)\\}}
  \7
\mu(A(x_1,x_2){\cap} A(y_1,y_2){\cap} S).
\end{equation}

Для доказательства это неравенства воспользуемся формулой
(\ref{pk}) и формулой включения-исключения
$$
\mu\left(\overline{M}_k{\cap}
S\right)=\mu\left(\bigcup_{\substack{
{(x_1,\,x_2)}{\in}\mathbb{Z}^2\\
1{\leqslant}\max (|x_1|,|x_2|){\leqslant}k\\
}} (A(x_1,x_2){\cap} S)\right)
{\geqslant}\mu\left(\bigcup\limits_{(x_1,\,x_2){\in} E_k}
(A(x_1,x_2){\cap} S)\right){\geqslant}$$
$$
{\geqslant}\sum\limits_{(x_1,\,x_2){\in} E_k}\mu(A(x_1,x_2){\cap}
S)\;-\7\5\sum_{\substack{
  (x_1,\,x_2){\in} E_k\\
  (y_1,\,y_2){\in} E_k\\
  (x_1,x_2)\neq(y_1,y_2)\\
  }}\5
\mu(A(x_1,x_2){\cap} A(y_1,y_2){\cap} S).
$$
Применив свойства множества $A(x_1,x_2)$, можно показать, что
выполняются равенства $$\mu(A(x_1,x_2){\cap}
S(x_1,x_2,k_1,k_2))=\left\{%
\begin{array}{ll}
\frac{2\varepsilon}{k^2}\frac{k_1k_2}{x_1x_2}, & \hbox{$x_1, x_2\neq 0$;} \\
\\
\frac{2\varepsilon}{k^2}\frac{k_2}{x_2}, & \hbox{$x_1=0, x_2\neq 0$;} \\
\\
\frac{2\varepsilon}{k^2}\frac{k_1}{x_1}, & \hbox{$x_2=0, x_1\neq 0$.} \\
\end{array}%
\right.$$

Оценим снизу первую сумму из неравенства (\ref{osn})
$$
\sum\limits_{(x_1,\,x_2){\in} E_k}\5\mu(A(x_1,x_2){\cap}
S)\geqslant\5\sum\limits_{(x_1,\,x_2){\in} E_k}\5\mu
(A(x_1,x_2){\cap}
S(x_1,x_2,k_1,k_2))=\frac{2\varepsilon}{k^2}\sum\limits_{(x_1,\,x_2){\in}
E_k}\frac{k_1k_2}{x_1x_2}\geqslant$$
$$
\geqslant\frac{2\varepsilon}{k^2} \sum\limits_{(x_1,\,x_2){\in}
E_k}\frac{(\lambda x_1-1)(\lambda x_2-1)}{x_1x_2}
=\frac{2\varepsilon}{k^2}\sum\limits_{(x_1,\,x_2){\in}
E_k}\left(\lambda^2-
\frac{\lambda}{x_1}-\frac{\lambda}{x_2}+\frac{1}{x_1x_2}\right)=
$$
$$
=\frac{2\varepsilon}{k^2}\left(\lambda^2\frac{k^2}{4\zeta(2)}+O(k\ln
k)\right)=\frac{\lambda^2\varepsilon}{2\zeta(2)}+O\left(\frac{\ln
k}{k}\right)\geqslant\frac{\lambda^2\varepsilon}{3\zeta(2)}.
$$

Оценим сверху меру множества $A(x_1,x_2){\cap} A(y_1,y_2){\cap}
S$. Квадрат $S$ и параллелограмм могут пересекаться, если
расстояние между их центрами меньше, чем
$\frac{\lambda}{\sqrt{2}}+\frac{2\varepsilon}{k}\leqslant\lambda$.
Если взять квадрат $3S$, у которого центр совпадает с центром
квадрата $S$ и сторона равна $3\lambda$, то центры тех
параллелограммов, которые задевают квадрат $S$ будут лежать внутри
квадрата $3S$. По лемме 1 будет
$$
\#(3S{\cap}\Lambda)\leqslant9\lambda^2\Delta+6\lambda\sqrt{x_1^2+y_1^2}+6\lambda\sqrt{x_2^2+y_2^2},
$$
где $\Delta=\left|\det\left(%
\begin{array}{cc}
  x_1 & x_2\\
  y_1 & y_2\\
\end{array}%
\right)\right|$.

Символ $\#(\cdot)$, здесь и далее, обозначает количество элементов
в множестве.

Для ограничения сверху второй суммы из неравенства (\ref{osn})
воспользуемся формулой (\ref{plr}) и леммой 2
$$
\sum_{\substack{
  (x_1,x_2){\in} E_k\\
  (y_1,y_2){\in} E_k\\
  (x_1,x_2)\neq(y_1,y_2)\\
  }}\7
\mu (A(x_1,x_2){\cap}
A(y_1,y_2){\cap}S)\leqslant\7\sum_{\substack{
  (x_1,x_2){\in} E_k\\
  (y_1,y_2){\in} E_k\\
  (x_1,x_2)\neq(y_1,y_2)\\
  }}\7\#(3S{\cap}\Lambda)\frac{4\varepsilon^2}{\Delta
k^4}\leqslant
$$
$$\leqslant\7\sum_{\substack{
  (x_1,x_2){\in} E_k\\
  (y_1,y_2){\in} E_k\\
  (x_1,x_2)\neq(y_1,y_2)\\
  }}\7
\left(9\lambda^2\Delta+6\lambda\sqrt{x_1^2+y_1^2}+6\lambda\sqrt{x_2^2+y_2^2}
\right)\frac{4\varepsilon^2}{\Delta k^4}\leqslant
$$
$$
\leqslant\frac{12\varepsilon^2}{k^4}\7\sum_{\substack{
  (x_1,x_2){\in}E_k\\
  (y_1,y_2){\in}E_k\\
  (x_1,x_2)\neq(y_1,y_2)\\
  }}\!\!\!
\left(3\lambda^2+\lambda\frac{2\sqrt{2}k}{\Delta}+\lambda\frac{2\sqrt{2}k}{\Delta}
\right)
=\frac{36\lambda^2\varepsilon^2}{k^4}\7\sum_{\substack{
  (x_1,x_2){\in}E_k\\
  (y_1,y_2){\in}E_k\\
  (x_1,x_2)\neq(y_1,y_2)\\
  }}\7 1+
\frac{48{\sqrt2}\lambda\varepsilon^2}{k^3}\7\sum_{\substack{
  (x_1,x_2){\in}E_k\\
  (y_1,y_2){\in}E_k\\
  (x_1,x_2)\neq(y_1,y_2)\\
  }}\7\frac{1}{\Delta}\leqslant
$$
$$
\leqslant\frac{36\lambda^2\varepsilon^2}{k^4}\7\sum_{\substack{
  (x_1,x_2){\in}E_k\\
  (y_1,y_2){\in}E_k\\
  }}\!\! 1+
\frac{432{\sqrt2}\lambda\varepsilon^2\ln
k}{k}\leqslant\frac{37\lambda^2\varepsilon^2}{\zeta^2(2)}.
$$
Ограничение снизу для $\mu(\overline{M}_k{\cap}S)$ доказано.

Обозначим $E^\ast_k=\{(x_1,x_2){\in}\mathbb{Z}^2:1{\leqslant}
x_1{\leqslant} k,1{\leqslant}|x_2|{\leqslant} k\}$ и
$S^\ast{=}S(x_1, x_2,k_1+1,k_2+1)$.

Воспользуемся формулой (\ref{pk}), свойством
$A(x_1,x_2){=}A({-}x_1,{-}x_2)$ и получим ограничение сверху на
$\mu(\overline{M}_k{\cap}S)$
$$
\mu(\overline{M}_k{\cap}
S)=\mu\left(\bigcup_{\substack{
{(x_1,x_2)}{\in}\mathbb{Z}^2\\
1{\leqslant}\max (|x_1|,|x_2|){\leqslant}k\\
}} (A(x_1,x_2){{\cap}}S)\right)
\leqslant{\mu}\left(\bigcup_{\substack{
                                  x_1=0\\
                                    1{\leqslant}x_2{\leqslant} k}}
\left(A(x
_1,x_2){{\cap}} S\right)\right)+
$$
$$
+\mu\left(\bigcup_{\substack{
x_2=0 \\
1{\leqslant}x_1{\leqslant} k}}
\left(A(x_1,x_2){{\cap}}S\right)\right)+\mu\left(\bigcup_{\substack{(x_1,\,x_2){\in}
E_k^\ast}}\left(A(x_1,x_2){{\cap}}S\right)\right){\leqslant}$$
$$
{\leqslant}2\sum_{\substack{
                                  x_1{=}0\\
                                  1{\leqslant}x_2{\leqslant}k \\
                                  }}
\mu\left(A(x_1, x_2){{\cap}}S^\ast\right){+}\sum\limits_{(x_1,\,
x_2){\in}E^\ast_k}\mu\left(A(x_1,x_2){{\cap}}S^\ast\right){=}
\frac{4\varepsilon}{k^2}\sum\limits_{x_2=1}^k \frac{k_2+1}{x_2}+
$$
$$
+\frac{2\varepsilon}{k^2}\sum\limits_{(x_1 ,\, x_2){\in}
E^\ast_k}\frac{(k_1+1)(k_2+1)}{x_1x_2}\leqslant
\frac{4\varepsilon}{k^2}\sum\limits_{x_2=1}^k \frac{\lambda
x_2+1}{x_2}+\frac{2\varepsilon}{k^2}\sum\limits_{(x_1,\, x_2){\in}
E^\ast_k}\frac{(\lambda x_1+1)(\lambda x_2+1)}{x_1x_2}=$$
$$
=\frac{4\varepsilon}{k^2}\sum\limits_{x_2=1}^k
\left(\lambda+\frac{1}{x_2}\right)+
\frac{2\varepsilon}{k^2}\sum\limits_{(x_1,\, x_2){\in}
E^\ast_k}\left(\lambda^2+
\frac{\lambda}{x_1}+\frac{\lambda}{x_2}+\frac{1}{x_1x_2}\right)=$$
$$
=\frac{2\varepsilon}{k^2}\left(2\lambda^2k^2+O(k\ln
k)\right)=4\varepsilon\lambda^2+O\left(\frac{\ln
k}{k}\right)\leqslant5\varepsilon\lambda^2.
$$
Лемма 3 доказана.

{\bf Лемма 4.} Для любого квадрата $S$ со стороной $\lambda$ верно
неравенство
$$
\lambda^2-5\lambda^2\varepsilon\leqslant\mu(S{\cap}\underline{M_k}).
$$
{\bf Доказательство.} Воспользуемся равенством (\ref{raven}) и
леммой 3
$$
\mu(S{\cap}\underline{M_k})=\mu(S)-\mu(S{\cap}\overline{M}_k)\geqslant\lambda^2-5\varepsilon\lambda^2.
$$
Лемма 4 доказана.

\textbf{3. Изучение функции $\psi_\Theta(t)$ в случае
$m\geqslant2$ и $n=1$.}

В этом случае функция $\psi_\Theta(t)$ имеет вид
$$
\psi_{\left(\scriptsize%
\begin{array}{c}
  \alpha_1\\
  \vdots\\
  \alpha_m\\
\end{array}%
\right)}(t)=\min\limits_{1\leqslant x\leqslant
t}\max\limits_{j=1,\ldots,m}\|x\alpha_j\|.
$$

На кубе $[0;1]^m$ рассмотрим два множества:
$$
\underline{M_k}=\left\{(\alpha_1,\ldots,\alpha_m){\in}[0;1]^m:
\psi_{\left(\scriptsize%
\begin{array}{c}
  \alpha_1\\
  \vdots\\
  \alpha_m\\
\end{array}%
\right)}(k)>\frac{\varepsilon}{\sqrt[m]{k}}\right\}=$$
$$
=\left\{(\alpha_1,\ldots,\alpha_m){\in}[0;1]^m:\forall q{\in}[1;k]
\max\limits_{1\leqslant i\leqslant m}\left\{\|\alpha_i
q\|\right\}>\frac{\varepsilon}{\sqrt[m]{k}}\right\}=$$
$$
=\left\{(\alpha_1,\ldots,\alpha_m){\in}[0;1]^m:\forall q{\in}[1;k]
\,\,\exists (p_1,\ldots,p_m){\in}\mathbb{Z}^m\,
\max\limits_{1\leqslant i\leqslant m}\left\{|\alpha_i
q-p_i|\right\}>\frac{\varepsilon}{\sqrt[m]{k}}\right\}=$$
$$
=\left\{(\alpha_1,\ldots,\alpha_m){\in}[0;1]^m:\forall
q{\in}[1;k]\,\,\exists (p_1,\ldots,p_m){\in}\mathbb{Z}^m\,
\max\limits_{1\leqslant i\leqslant
m}\left\{\left|\alpha_i-\frac{p_i}{q}\right|\right\}>
\frac{\varepsilon}{q\sqrt[m]{k}}\right\};
$$
$$
\overline{M}_k=\left\{(\alpha_1,\ldots,\alpha_m){\in}[0;1]^m:
\psi_{\left(\scriptsize%
\begin{array}{c}
  \alpha_1\\
  \vdots\\
  \alpha_m\\
\end{array}%
\right)}(k)\leqslant\frac{\varepsilon}{\sqrt[m]{k}}\right\}=$$
$$
=\left\{(\alpha_1,\ldots,\alpha_m){\in}[0;1]^m:\exists q{\in}[1;k]
\max\limits_{1\leqslant i\leqslant m}\left\{\|\alpha_i
q\|\right\}\leqslant\frac{\varepsilon}{\sqrt[m]{k}}\right\}=$$
$$
=\left\{(\alpha_1,\ldots,\alpha_m){\in}[0;1]^m:\exists q{\in}[1;k]
\,\, \forall (p_1,\ldots,p_m){\in}\mathbb{Z}^m
\max\limits_{1\leqslant i\leqslant
m}\left\{|\alpha_iq-p_i|\right\}\leqslant\frac{\varepsilon}{\sqrt[m]{k}}\right\}=$$
$$
=\left\{(\alpha_1,\ldots,\alpha_m){\in}[0;1]^m:\exists q{\in}[1;k]
\,\, \forall(p_1,\ldots,p_m){\in}\mathbb{Z}^m
\max\limits_{1\leqslant i\leqslant
m}\left\{\left|\alpha_i-\frac{p_i}{q}\right|\right\}
\leqslant\frac{\varepsilon}{q\sqrt[m]{k}}\right\}.
$$

Верно равенство
\begin{equation}\label{vran}
\mu(\underline{M_k})=1-\mu(\overline{M}_k).
\end{equation}

Для числа $q$ определим множество $A(q)$, как объединение кубов с
центрами в точках
$\left(\frac{p_1}{q},\frac{p_2}{q},\ldots,\frac{p_m}{q}\right)$
при $p_i=0,1,\ldots,q$ и  со стороной
$\frac{2\varepsilon}{q\sqrt[m]{k}}$.

Верно равенство
\begin{equation}\label{tr}
\overline{M}_k=\bigcup\limits_{q=1}^k
\left(A(q){\cap}[0,1]^m\right).
\end{equation}

Символ $\lceil \cdot\rceil$ обозначает целую часть сверху.

{\bf Лемма 5.} Для любого куба $S$ со стороной $\lambda$
выполняется неравенство
$$
\lambda^m-2(4\lambda\varepsilon)^m\leqslant\mu(S{\cap}
\underline{M_k}).
$$

{\bf Доказательство.} Проекция множества $A(q)$ на каждую
координатную ось это объединение отрезков вида
$\left[\frac{p}{q}-\frac{\varepsilon}{q\sqrt[m]{k}},\frac{p}{q}+\frac{\varepsilon}{q\sqrt[m]{k}}\right]$
при $p=0, 1, \ldots,  q$. Очевидно отрезок длины $\lambda$ заденет
не более
$\left\lceil\left(\lambda+\frac{2\varepsilon}{q\sqrt[m]{k}}\right):\frac1q\right\rceil\leqslant
\left(\lambda+\frac{2\varepsilon}{q\sqrt[m]{k}}\right)q+1\leqslant\lambda
q+2$ отрезков из проекции множества $A(q)$ по каждой оси. Получаем
неравенство
$$
\mu(S{\cap} A(q))\leqslant(\lambda
q+2)^m\frac{(2\varepsilon)^m}{kq^m}.
$$

Применим формулу (\ref{tr}) и найдем ограничение сверху на меру
пересечения
\begin{gather*}
\mu(S{\cap}\overline{M}_k)=\mu\left(\bigcup\limits_{q=1}^k(S{\cap}
A(q))\right)\leqslant\sum\limits_{q=1}^k\mu(S{\cap}
A(q))\leqslant\sum\limits_{q=1}^k\left(\lambda
q+2\right)^m\frac{(2\varepsilon)^m}{kq^m}\leqslant\\
\leqslant\sum\limits_{q=1}^{k}\left((2\lambda
q)^m+4^m\right)\frac{(2\varepsilon)^m}{kq^m}=\frac{(4\lambda\varepsilon)^m}{k}\sum\limits_{q=1}^k1
+\frac{8^m\varepsilon^m}{k}\sum\limits_{q=1}^{k}
\frac{1}{q^m}\leqslant\\
\leqslant(4\lambda\varepsilon)^m+
\frac{8^m\varepsilon^m}{k}\sum\limits_{q=1}^{\infty}
\frac{1}{q^2}\leqslant2(4\lambda\varepsilon)^m.
\end{gather*}
Применим полученное неравенство и формулу (\ref{vran})
$$
\mu(S{\cap}\underline{M_k})=\mu(S)-\mu(S{\cap}\overline{M}_k)\geqslant\lambda^m-2(4\lambda\varepsilon)^m.
$$
Лемма 5 доказана.

{\bf Лемма 6.} Для $m\geqslant2$ верно неравенство
$$
\sum\limits_{\lceil k/2\rceil\leqslant q_1<q_2\leqslant
k}\!\frac{(q_1,q_2)^m}{q_2^m}\leqslant \frac{2k}{5}.
$$

{\bf Доказательство.} Пусть $(q_1,q_2)=d$, $q_1{=}dl$
$\left(l{\geqslant}\left\lceil\frac{k}{2d}\right\rceil\right)$  и
$q_2{=}dp$ $\left(p{\geqslant}
l{+}1{\geqslant}\left\lceil\frac{k}{2d}\right\rceil{+}1\right)$ с
условием $(l,p)=1$, тогда
$$
\sum\limits_{\lceil k/2\rceil\leqslant q_1<q_2\leqslant
k}\frac{(q_1,q_2)^ m}{q_2^m}=\sum\limits_{q_2=\left\lceil\frac
k2\right\rceil{+}1}^k\sum\limits_{q_1=\left\lceil\frac
k2\right\rceil}^{q_2-1}\frac{(q_1,q_2)^
m}{q_2^m}=\sum\limits_{d=1}^{[k/2]}
\sum\limits_{p=\left\lceil\frac{k}{2d}\right\rceil{+}1}^{\left[
k/d\right]}\sum\limits_{\scriptsize{\begin{array}{c}
l=\left\lceil\frac{k}{2d}\right\rceil\\
(l,p){=}1\\
\end{array}}}^{p-1}
\frac{1}{p^m}\leqslant$$
$$
{\leqslant}\sum\limits_{d=1}^{[k/2]}\sum\limits_{\scriptsize
p=\left\lceil\frac{k}{2d}\right\rceil{+1}}^{\left[
k/d\right]}\sum\limits_{\scriptsize{\begin{array}{c}
l{=}1\\
(l,p){=}1\\
\end{array}}}^{p-1}
\frac{1}{p^m}{=}\sum\limits_{d=1}^{[k/2]} \sum\limits_{\scriptsize
p=\left\lceil\frac{k}{2d}\right\rceil+1}^{[k/d]}\frac{\varphi(p)}{p^
m}\leqslant\sum\limits_{p=2}^k\sum\limits_{d=1}^{[k/p]}\frac{\varphi(p)}{p^
m}=
$$
$$
=\sum\limits_{p=2}^k\left[\frac{k}{p}\right]\!\!\frac{\varphi(p)}{p^
m}{\leqslant}k\!\sum\limits_{p=2}^k\frac{\varphi(p)}{p^{m+1}
}{\leqslant}
k\left(\sum\limits_{p=1}^\infty\frac{\varphi(p)}{p^{m+1}}
{-}1\right){=}k\left(\frac{\zeta(m)}{\zeta(m+1)}-1\right){\leqslant}
k\left(\frac{\zeta(2)}{\zeta(3)}{-}1\right)\!{\leqslant}\frac{2k}{5}.
$$

Доказательство равенства
$$
\sum\limits_{p=1}^\infty\frac{\varphi(p)}{p^{m}}
=\frac{\zeta(m-1)}{\zeta(m)}
$$
можно найти в \cite{Hardy} на стр. 250.

Лемма 6 доказана.

{\bf Лемма 7.} Пусть $W$ - выпуклая область на плоскости, $N$ -
число целых точек в области $W$. Если найдутся 3 целые точки
внутри области не лежащие на одной прямой, то
$$
N\leqslant 2\mu(W)+2.
$$
{\bf Доказательство.} Рассмотрим все целые точки, лежащие в
области $W$, для них построим выпуклую оболочку. Получим выпуклый
многоугольник. Применив формулу Пика, получим неравенство
$$
B+\frac{\Gamma}{2}-1\leqslant\mu(W).
$$
где $B$ есть количество целочисленных точек внутри многоугольника,
а $\Gamma$ — количество точек на границе многоугольника. Из
соотношения $N=B+\Gamma$ получаем
$$
N-2\leqslant2\mu(W).
$$
Лемма 7 доказана.

Рассмотрим плоскость $Op_1p_2$. Для пары
$(q_1,q_2){\in}\mathbb{N}^2$ и числа $\delta>0$ введем
обозначения:

$C$ --- прямоугольник со сторонами $(1+\delta)\lambda q_1$ и
$(1+\delta)\lambda q_2$, параллельными осям;\\

$C_0=\left[0,\frac{(1+\delta)}{2}\lambda
q_1\right]\times\left[0,\frac{(1+\delta)}{2}\lambda q_2\right]$;\\

$D=\left\{(p_1,p_2)\in\mathbb{R}^2\colon
\left|q_2p_1-q_1p_2\right|\leqslant
d\left[\frac{\varepsilon(q_1+q_2)}{d\sqrt[m]{k}}\right]\right\}$,
где $d$ --- наибольший общий делитель чисел $q_1$ и $q_2$.\\

$l_q$ --- прямая $q_2p_1-q_1p_2=q$.

Запишем несколько свойств прямых вида $l_q$:

1. Расстояние между двумя соседними прямыми вида $l_q$ равно
$\frac{1}{\sqrt{q_1^2+q_2^2}}$.

2. Все целые точки $(p_1,p_2){\in}\mathbb{Z}^2$ лежат на прямых
вида $l_q$, расстояние между двумя соседними точками на одной
прямой равно $\sqrt{q_1^2+q_2^2}$.\\

Запишем некоторые свойства области $D$:

1. Количество прямых $l_q$ в области $D$ равно $
2d\left[\frac{\varepsilon\left(q_1+q_2\right)}{d\sqrt[m]{k}}\right]+1;
$

2. Ширина полосы $D$ равна $
h=2\left[\frac{\varepsilon\left(q_1+q_2\right)}{d\sqrt[m]{k}}\right]\frac{d}{\sqrt{q_1^2+q_2^2}}\leqslant
\frac{2\sqrt{2}\varepsilon}{\sqrt[m]{k}}. $

3. Если $d>2\varepsilon k^{1-1/m}$, то $D=l_0$.\\

Определим множества:

$Tr_k=\{(q_1,q_2){\in}\mathbb{Z}^2\colon\, \lceil k/2\rceil\leqslant
q_1<q_2\leqslant
k\}$;\\

$J=\{(q_1,q_2){\in} Tr_k\colon\exists C=C_1\,\, \#(D{\cap}
C_1{\cap}\mathbb{Z}^2)\geqslant N_0\}$, где $N_0=8(1+\delta)\varepsilon\lambda k^{1-1/m}+6;$\\

$J_0=\{(q_1,q_2){\in} Tr_k\colon (q_1,q_2)=d>d_0\}$, где
$d_0=9\varepsilon k^{1-1/m}$;\\

$J_1=\left\{(q_1,q_2){\in}Tr_k\colon
\#(D{\cap}C_0{\cap}\mathbb{Z}^2)\geqslant\frac{N_0}{2}\right\}$;\\

$Pr_k{=}\left\{(p_1,p_2){\in}
\left[0,\frac{\sqrt[m]{k}}{4\varepsilon}\right]^2\!\cap\mathbb{N}^2\colon\,(p_1,p_2)=1,
1\leqslant\frac{p_2}{p_1}\leqslant2\right\}.$

{\bf  Лемма 8.} Для каждой пары $(q_1,q_2){\in} J$ все целые точки
в $D{\cap} C_1$ лежат на одной прямой $l$ c шагом $\gamma$, где
$$
1\leqslant\gamma\leqslant\frac{\sqrt[m]{k}}{4\varepsilon}.
$$

{\bf Доказательство.} Область $D{\cap} C_1$ можно вложить в
параллелограмм с высотой $h$ и стороной
$\sqrt{2}(1{+}\delta)\lambda k$. Таким образом мера
$\mu(D{{\cap}}C_1)\leqslant4(1{+}\delta)\varepsilon\lambda
k^{1{-}1/m}$.

Пусть точки не лежат на одной прямой, тогда по лемме 7, количество
целых точек $N$ в области $\mu(D{\cap} C_1)$ будет
$$
N\leqslant2\mu(D{\cap}
C_1)+2\leqslant8(1+\delta)\varepsilon\lambda k^{1{-}1/m}+2,
$$
получаем противоречие с тем, что $(q_1,q_2){\in} J$.

Максимальный отрезок в $D{\cap} C_1$ это диагональ прямоугольника
$C_1$. Длина диагонали равна $\sqrt{2}(1+\delta)\lambda k$,
минимальное количество точек на отрезке такой длины должно быть не
менее $N_0$, поэтому шаг $\gamma$ не может превышать величины
$$
\frac{\sqrt{2}(1+\delta)\lambda
k}{N_0-1}\leqslant\frac{\sqrt{2}(1+\delta)\lambda
k}{N_0-6}=\frac{\sqrt{2}\sqrt[m]{k}}{8\varepsilon}\leqslant\frac{\sqrt[m]{k}}{4\varepsilon}.
$$
Лемма 8 доказана.

{\bf Лемма 9.} Пусть $l$ - прямая из леммы 8, $\gamma$ -
расстояние между двумя соседними целыми точками на этой прямой,
$\omega$ - угол между прямыми $l$ и $l_0$, тогда
$$
\sin\omega\leqslant\frac{1}{2k\lambda\gamma}.
$$

{\bf Доказательство.} Определим число
$a=d\left[\frac{\varepsilon}{d\sqrt[m]{k}}\left(q_1+q_2\right)\right]$.
Построим прямоугольный треугольник $A_1A_2A_3$ с гипотенузой
$A_1A_2$, где

$A_1=l{\cap} l_a$;

$A_2=l{\cap} l_{-a}$;

$A_3{\in} l_{a}$;

$A_1A_3\bot l$.

Длины сторон $|A_1A_2|\geqslant (N_0-1)\gamma$ и $|A_1A_3|=h$.
Оценим синус угла между прямыми $l$ и $l_0$
$$
\sin\omega=\frac{|A_1A_3|}{|A_1A_2|}\leqslant\frac{h}{(N_0-1)\gamma}
\leqslant\frac{2\sqrt{2}\varepsilon}{\sqrt[m]{k}(N_0-6)\gamma}\leqslant
\frac{1}{2k\lambda\gamma}.
$$
Если $D=l_0$, то $l=l_0$ и $\sin\omega=0$.

Лемма 9 доказана.

{\bf Лемма 10.} Если $(q_1,q_2){\in}J$,  то
$$
\#(D{\cap} C_0{\cap}\mathbb{Z}^2)\geqslant\frac{N_0}{2}.
$$

{\bf Доказательство.} Пусть $(q_1,q_2){\in} J$, тогда существует
квадрат $C_1$, такой что в области $D{\cap} C_1$ будет лежать
$N\geqslant N_0$ целых точек. Из всех целых точек в $D{\cap}C_1$
выделим две точки $p_s$ и $p_d$ - самая левая и самая правая
соответственно по оси $Op_1$. Возьмем точку
$p_c=\frac{p_s+p_d}{2}$. Рассмотрим 2 случая.

1-ый случай. Если $N$ - нечетное, то $p_c$ - целая точка.
Обозначим через $\phi$ сдвиг на вектор $-\mathbf{p}_c$, где вектор
$\mathbf{p}_c$ - радиус-вектор точки $p_c$.

Точки $\phi(p_c)=(0,0){\in}D{\cap}C_0$ и
$\phi(p_d){\in}D{\cap}C_0$, а значит и образы всех целых точек
между $p_c$ и $p_d$ тоже лежат в области $D{\cap}C_0$. Количество
целых точек в области $D{\cap}C_0$ будет
$$
\#(D{\cap}C_0{\cap}\mathbb{Z}^2)\geqslant
\frac{N+1}{2}\geqslant\frac{N_0+1}{2}>\frac{N_0}{2}.
$$

2-ой случай. Если $N$ - четное, то точка $p_c$ лежит посередине
двух целых точек с прямой $l$, возьмем правую точку, обозначим её
$p_{cd}$. И сделаем сдвиг на вектор $-\mathbf{p}_{cd}$, где вектор
$\mathbf{p}_{cd}$ - радиус-вектор точки $p_{cd}$. Точки
$\phi(p_{cd})=(0,0){\in}D{\cap}C_0$ и
$\phi(p_{d}){\in}D{\cap}C_0$, а значит и образы всех целых точек
между $p_{cd}$ и $p_d$ тоже лежат в области $D{\cap}C_0$.

Количество целых точек в области $D{\cap}C_0$ будет
$$
\#(D{\cap}C_0{\cap}\mathbb{Z}^2)\geqslant
\frac{N}{2}\geqslant\frac{N_0}{2}.
$$

Лемма 10 доказана.

{\bf Лемма 11.} Для любого $(q_1,q_2){\in} Tr_k$ выполняется
неравенство
$$
\#(D{\cap} C{\cap}\mathbb{Z}^2)\leqslant\left(\left(32\varepsilon
\lambda k^{1-1/m}\right)^m+ \left(\frac{16\varepsilon
k^{1-1/m}}{d}\right)^m+(8\lambda d)^m+4^m\right)^{1/m}.
$$
{\bf Доказательство.}

Количество прямых $l_q$ в области $D$, на которых лежат целые
точки равно
$$
2\left[\frac{\varepsilon}{d\sqrt[m]{k}}\left(q_1+q_2\right)\right]+1\leqslant\frac{4\varepsilon
k^{1-1/m}}{d}+1.
$$

Количество точек на каждой прямой не превосходит количества точек
на диагонали прямоугольника $C$. На каждой прямой расстояние между
точками $\frac{\sqrt{q_1^2+q_2^2}}{d}$. Длина диагонали
$\sqrt{\left((1+\delta)\lambda
q_1\right)^2+\left((1+\delta)\lambda
q_2\right)^2}=(1+\delta)\lambda\sqrt{q_1^2+q_2^2}$. Значит
наибольшее количество точек на каждой прямой не превосходит
величины
\begin{equation}\label{line}
\frac{(1+\delta)\lambda\sqrt{q_1^2+q_2^2}}{\frac{\sqrt{q_1^2+q_2^2}}{d}}+1=(1+\delta)\lambda
d+1.
\end{equation}

Таким образом всё количество целых точек в области $D{\cap}C$
оценим так
$$
\#(D{\cap} C{\cap}\mathbb{Z}^2)\leqslant\left(\frac{4\varepsilon
k^{1-1/m}}{d}+1\right)\left(2\lambda
d+1\right)=\left(\left(\frac{4\varepsilon k^{1-1/m}}{d}+1\right)^m
\left(2\lambda d+1\right)^m\right)^{1/m}\leqslant
$$
$$
\leqslant \left(\left(\left(\frac{8\varepsilon
k^{1-1/m}}{d}\right)^m+2^m\right) \left((4\lambda
d)^m+2^m\right)\right)^{1/m}=
$$
$$
=\left(\left(32\varepsilon \lambda k^{1-1/m}\right)^m+
\left(\frac{16\varepsilon k^{1-1/m}}{d}\right)^m+(8\lambda
d)^m+4^m\right)^{1/m}.
$$
Лемма 11 доказана.

Вокруг каждой точки $(p_1,p_2){\in} Pr_k$ построим сектор
$\mathrm{V}(p_1,p_2)$ с углом $\phi$, для которого верно
неравенство $\sin\phi=\frac{1}{k\lambda\gamma}$, где
$\gamma=\sqrt{p_1^2+p_2^2}$ и точка $(p_1,p_2)$ лежит на
биссектрисе угла $\phi$. Под сектором будем понимать область между
двумя прямыми вида $y=k_1x$ и $y=k_2x$. Объединение всех
получившихся секторов обозначим через~$\mathrm{V}$.

{\bf Лемма 12.} Верны включения $J_0\subset J\subset
J_1\subset\left(\mathrm{V}{\cap}Tr_k\right).$

{\bf Доказательство.} Докажем включение $J_0\subset J$. Если пара
$(q_1,q_2){\in}J_0$, то $D=l_0$. Возьмем прямоугольник
$C=\left[0,(1+\delta)\lambda
q_1\right]\times\left[0,(1+\delta)\lambda q_2\right]$. Оценим
количество целых точек в $D{\cap}C$
$$
\#\left(D{\cap}C{\cap}\mathbb{Z}^2\right)=\#(l_0{\cap}C{\cap}\mathbb{Z}^2)=
\left[\frac{(1+\delta)\lambda\sqrt{q_1^2+q_2^2}}
{\frac{\sqrt{q_1^2+q_2^2}}{d}}\right]+1>\frac{(1+\delta)\lambda\sqrt{q_1^2+q_2^2}}
{\frac{\sqrt{q_1^2+q_2^2}}{d}}=$$
$$
=(1+\delta)\lambda d>9(1+\delta)\varepsilon\lambda
k^{1-1/m}>8(1+\delta)\varepsilon\lambda k^{1-1/m}+6.
$$
Включение $J_0\subset J$ доказано.

Включение $J\subset J_1$ доказано в лемме 10.

Если $(q_1,q_2){\in}J_1$, то в области $D{\cap}C_0$ все целые
точки лежат на одной прямой $l$, так как
$\mu(D{\cap}C_0)\leqslant2(1{+}\delta)\varepsilon\lambda
k^{1{-}1/m}$.

Пусть точки не лежат на одной прямой, тогда по лемме 7, количество
целых точек в области $\mu(D{\cap} C_0)$ будет
$$
2\mu(D{\cap} C_0)+2\leqslant4(1+\delta)\varepsilon\lambda
k^{1{-}1/m}+2<\frac{N_0}{2},
$$
получаем противоречие с тем, что $(q_1,q_2){\in}J_1$. Прямые $l$ и
$l_0$ пересекаются в начале координат и синус угла между прямыми
$l$ и $l_0$ не превосходит величины $\frac{1}{2k\lambda\gamma}$
Для доказательства надо повторить рассуждения, аналогичные
рассуждениям в лемме 9.

Докажем включение $J_1\subset\left(\mathrm{V}{\cap}Tr_k\right)$.

Если точка $(q_1,q_2){\notin}\mathrm{V}$, то
$(q_1,q_2){\notin}J_1$, потому как синус угла между прямой $l_0$ и
любым радиус-вектором $(p_1,p_2)$, где $(p_1,p_2){\in}Pr_k$ больше
чем
$\sin\frac{\phi}{2}>\frac{\sin\phi}{2}=\frac{1}{2k\lambda\gamma}$,
что противоречит тому что синус угла между прямыми $l$ и $l_0$ не
превосходит величины $\frac{1}{2k\lambda\gamma}$.

Лемма 12 доказана.

{\bf Лемма 13.} Количество пар $(q_1,q_2){\in}J$ не более
$\frac{2k^{1+\frac{3}{2m}}}{\lambda\varepsilon^2}$.

{\bf Доказательство.} Оценим сколько точек $(q_1,q_2){\in}Tr_k$
попадают в каждый сектор $\mathrm{V}(p_1,p_2)$. Возможны 2 случая.

Первый случай. Все целые точки лежат на прямой $l_0$, тогда их
количество не больше чем $\frac{\sqrt{2}k}{\gamma}+1$.

Второй случай. Для области $\mathrm{V}(p_1,p_2){\cap}[1,k]^2$
выполняются условия леммы 7, тогда количество точек не больше чем
$$
2\mu(\mathrm{V}(p_1,p_2){\cap}[1,k]^2){+}2\leqslant2\cdot\frac{1}{2}\left(\sqrt{2}k\right)^2\sin\phi+2
\leqslant\frac{2k}{\lambda\gamma}+2.
$$
Оценим общее количество точек
$$
\#J{\leqslant}\#(\mathrm{V}{\cap}
Tr_k){=}\#\left(\bigcup\limits_{(p_1,p_2){\in} Pr_k}\5
(\mathrm{V}(p_1,p_2){\cap} Tr_k)\right){\leqslant}
$$
$$
\leqslant\sum\limits_{(p_1,p_2){\in}
Pr_k}\#(\mathrm{V}(p_1,p_2){\cap}
Tr_k)\leqslant\sum\limits_{(p_1,p_2){\in}
Pr_k}\left(\frac{\sqrt{2}k}{\gamma}+1+\frac{2k}{\lambda\gamma}+2\right)\leqslant
$$
$$
\leqslant\sum\limits_{(p_1,p_2){\in}
Pr_k}\left(\frac{2k}{\lambda\gamma}+\frac{2k}{\lambda\gamma}+3\right)\leqslant\sum\limits_{(p_1,p_2){\in}
Pr_k}\frac{4k}{\lambda\gamma}+\frac{3k^{2/m}}{16\varepsilon^2}=
$$
$$
=\5\sum_{\substack{
(p_1,p_2){\in}Pr\\
1{\leqslant}\gamma{\leqslant}\frac{\sqrt[2m]{k}}{4\varepsilon}\\
}}\frac{4k}{\lambda\gamma}+\5\sum_{\substack{(p_1,p_2){\in}
Pr_k\\
\frac{\sqrt[2m]{k}}{4\varepsilon}
{<}\gamma{\leqslant}\frac{\sqrt[m]{k}}{4\varepsilon}\\
}}\frac{4k}{\lambda\gamma}+\frac{3k^{2/m}}{16\varepsilon^2}
\leqslant\5\sum_{\substack{(p_1,p_2){\in}
Pr\\
1\leqslant\gamma\leqslant\frac{\sqrt[2m]{k}}{4\varepsilon}\\
}}\5\frac{4k}{\lambda}+\5\sum_{\substack{(p_1,p_2){\in}
Pr\\
\frac{\sqrt[2m]{k}}{4\varepsilon}<\gamma\leqslant\frac{\sqrt[m]{k}}{4\varepsilon}\\
}}\frac{16k\varepsilon}{\lambda\sqrt[2m]{k}}+\frac{3k^{2/m}}{16\varepsilon^2}\leqslant
$$
$$
\leqslant\frac{4k}{\lambda}\frac{k^{1/m}}{16\varepsilon^2}+\frac{16k\varepsilon}{\lambda\sqrt[2m]{k}}
\frac{k^{2/m}}{16\varepsilon^2}+\frac{3k^{2/m}}{16\varepsilon^2}=\frac{k^{1+\frac1m}}{4\lambda\varepsilon^2}+
\frac{k^{1+\frac{3}{2m}}}{\lambda\varepsilon}+\frac{3k^{\frac2m}}{16\varepsilon^2}
\leqslant\frac{2k^{1+\frac{3}{2m}}}{\lambda\varepsilon^2}.
$$

Лемма 13 доказана.

{\bf Лемма 14.} Для любого куба $S$ со стороной $\lambda$
выполняется неравенство
$$
\frac{(2\lambda\varepsilon)^m}{6}-\frac{(68\lambda\varepsilon^2)^m}{4}\leqslant\mu(S{\cap}
\overline{M}_k).
$$
{\bf Доказательство.}

Ограничение снизу на меру пересечения множества $\overline{M}_k$ с
кубом S записывается неравенством
\begin{equation}\label{ff}
\mu(S{\cap}\overline{M}_k)\geqslant\sum\limits_{q=\lceil
k/2\rceil}^k\mu\left(S{\cap}
A(q)\right)-\7\sum\limits_{\small{\lceil k/2\rceil\leqslant
q_1<q_2\leqslant k}} \mu\left(S{\cap} A(q_1){\cap} A(q_2)\right).
\end{equation}
Применим формулу (\ref{tr}) и докажем это неравенство
$$
\mu(S{\cap}\overline{M}_k)=\mu\left(\bigcup\limits_{q=1}^{k}(S{\cap}
A(q))\right)\geqslant\mu\left(\bigcup\limits_{q=\lceil
k/2\rceil}^{k}(S{\cap} A(q))\right)\geqslant
$$
$$
\geqslant\sum\limits_{q=\lceil k/2\rceil}^k\mu\left(S{\cap}
A(q)\right)-\sum\limits_{\small{\lceil k/2\rceil\leqslant
q_1<q_2\leqslant k}} \mu\left((S{\cap} A(q_1))\bigcap (S{\cap}
A(q_2))\right)=
$$
$$
=\sum\limits_{q=\lceil k/2\rceil}^k\mu\left(S{\cap}
A(q)\right)-\sum\limits_{\small{\lceil k/2\rceil\leqslant
q_1<q_2\leqslant k}} \mu\left(S{\cap} A(q_1){\cap} A(q_2)\right).
$$

Первую сумму из (\ref{ff}) снизу можно оценить таким образом
\begin{equation}\label{ch1}
\sum\limits_{q=\lceil k/2\rceil}^k\mu(S{\cap} A(q))
\geqslant\frac{(1-\delta)^m(2\lambda\varepsilon)^m}{2}.
\end{equation}

Докажем это неравенство. Проекция множества $A(q)$ на каждую
координатную ось это объединение отрезков вида
$\left[\frac{p}{q}-\frac{\varepsilon}{q\sqrt[m]{k}},\frac{p}{q}+\frac{\varepsilon}{q\sqrt[m]{k}}\right]$
при $p=0, 1, \ldots,  q$. Минимальное количество проекций
множества $A(q)$, которые полностью попадут в отрезок длины
$\lambda$ равно
$\left[\left(\lambda-\frac{2\varepsilon}{\sqrt[m]{k}q}\right):\frac{1}{q}\right]\geqslant
\lambda q-\frac{2\varepsilon}{\sqrt[m]{k}}-1>\lambda q-2$.
Получаем неравенство
$$
\sum\limits_{q=\lceil k/2\rceil}^k\mu(S{\cap} A(q))\geqslant
\sum\limits_{q=\lceil k/2\rceil}^k\left(\lambda
q-2\right)^m\frac{2^m\varepsilon^m}{kq^m}
\geqslant\sum\limits_{q=\lceil k/2\rceil}^k\left((1-\delta)\lambda
q\right)^m\frac{2^m\varepsilon^m}{kq^m}=
$$
$$
=\sum\limits_{q=\lceil
k/2\rceil}^k\frac{(1-\delta)^m(2\lambda\varepsilon)^m}{k}\geqslant\frac{(1-\delta)^m(2\lambda\varepsilon)^m}{2}.
$$
В сумме $\sum\limits_{\small{\lceil k/2\rceil\leqslant
q_1<q_2\leqslant k}} \mu\left(S{\cap} A(q_1)\bigcap A(q_2)\right)$
найдем ограничение сверху на каждое слагаемое. Сначала найдем
ограничение на количество пересечений. Выпишем условия на попарное
пересечение множеств. Рассмотрим проекцию на первую ось.

Проекции множеств $A(q_1)$ и $S$ пересекаются, если
$$
\alpha_1-\frac{\varepsilon}{\sqrt[m]{k}q_1}\leqslant\frac{p_1}{q_1}<
\alpha_1+\lambda+\frac{\varepsilon}{\sqrt[m]{k}q_1}.
$$

Проекции множеств $A(q_2)$ и $S$ пересекаются, если
$$
\alpha_1-\frac{\varepsilon}{\sqrt[m]{k}q_2}\leqslant\frac{p_2}{q_2}<
\alpha_1+\lambda+\frac{\varepsilon}{\sqrt[m]{k}q_2}.
$$

Проекции множеств $A(q_1)$ и $A(q_2)$ пересекаются, если
существуют $p_1$ и $p_2$, такие что
$$
\left|\frac{p_1}{q_1}-\frac{p_2}{q_2}\right|
<\frac{\varepsilon}{\sqrt[m]{k}}\left(\frac{1}{q_1}+\frac{1}{q_2}\right).
$$
Если пара $(p_1,p_2)$ удовлетворяет всем трем неравенствам, то
проекции множеств $S$, $A(q_1)$ и $A(q_2)$ пересекаются. Перепишем
эти неравенства по другому
$$\left\{%
\begin{array}{ll}
    0\leqslant p_1-\alpha_1q_1+\frac{\varepsilon}{\sqrt[m]{k}}<\lambda
q_1+\frac{2\varepsilon}{\sqrt[m]{k}}\\
\\
    0\leqslant p_2-\alpha_1q_2+\frac{\varepsilon}{\sqrt[m]{k}}<\lambda
q_2+\frac{2\varepsilon}{\sqrt[m]{k}}\\
\\
    \left|p_1q_2-p_2q_1\right|<\frac{\varepsilon}{\sqrt[m]{k}}\left(q_1+q_2\right)\\
\end{array}%
\right.
$$

Рассмотрим плоскость $Op_1p_2$. Первые два неравенства задают
прямоугольник со сторонами $\lambda
q_1+\frac{2\varepsilon}{\sqrt[m]{k}}$ и $\lambda
q_2+\frac{2\varepsilon}{\sqrt[m]{k}}$, который содержится в
прямоугольнике $C$. Третье неравенство задает прямые вида
$p_1q_2{-}p_2q_1{=}q$, где $q=0,\pm d,\ldots, \pm
d\left[\frac{\varepsilon}{d\sqrt[m]{k}}\left(q_1{+}q_2\right)\right]$,
то есть в точности прямые из области $D$. Отсюда следует, что
количество решений системы будет не более, чем $\#(D{\cap}
C{\cap}\mathbb{Z}^2)$.

Всё множество пар $(q_1,q_2)$ разобьем на три случая:

пара $(q_1,q_2){\in} J_0$, тогда $D=l_0$ и по формуле (\ref{line})
$$
\#(D{\cap} C{\cap}\mathbb{Z}^2)\leqslant(1+\delta)\lambda d+1;
$$

пара $(q_1,q_2){\in} J\setminus J_0$, воспользуемся оценкой из
леммы 11
$$
\#(D{\cap} C{\cap}\mathbb{Z}^2)\leqslant\left(\left(32\varepsilon
\lambda k^{1-1/m}\right)^m+ \left(\frac{16\varepsilon
k^{1-1/m}}{d}\right)^m+(8\lambda d)^m+4^m\right)^{1/m};
$$

пара $(q_1,q_2){\in} Tr_k\setminus J$, тогда
$$
\#(D{\cap} C{\cap}\mathbb{Z}^2)<N_0=8(1+\delta)\varepsilon\lambda
k^{1-1/m}+6\leqslant17\varepsilon\lambda k^{1-1/m}.
$$
Сумму можно записать
$$
\sum\limits_{\lceil k/2\rceil\leqslant q_1<q_2\leqslant
k}\mu\left(S{\cap} A(q_1){\cap}
A(q_2)\right)=\left(\sum\limits_{J_0}+\sum\limits_{J\setminus
J_0}+\sum\limits_{Tr_k\setminus J}\right)\mu\left(S{\cap}
A(q_1){\cap} A(q_2)\right).
$$

Для удобства вместо $\sum\limits_{(q_1,q_2){\in}J}$,
$\sum\limits_{(q_1,q_2){\in}J\setminus J_0}$ и
$\sum\limits_{(q_1,q_2){\in}Tr_k}$ будем писать $\sum\limits_{J}$,
$\sum\limits_{J\setminus J_0}$ и $\sum\limits_{Tr_k}$
соответственно.

Применим лемму 6 и получим ограничение на сумму по
$(q_1,q_2){\in}J_0$
$$
\sum\limits_{J_0}\mu\left(S{\cap} A(q_1){\cap}
A(q_2)\right)\leqslant\sum\limits_{J_0}\left(\#(D{\cap}
C{\cap}\mathbb{Z}^2)\right)^m\frac{(2\varepsilon)^m}{kq_2^m}{\leqslant}
\sum\limits_{J_0}\left((1{+}\delta)\lambda
d{+}1\right)^m\frac{(2\varepsilon)^m}{kq_2^m}\leqslant
$$
$$
\leqslant\sum\limits_{J_0}\left((1{+}2\delta)\lambda
d\right)^m\frac{(2\varepsilon)^m}{kq_2^m}\leqslant
\frac{(1+2\delta)^m\lambda^m(2\varepsilon)^m}{k}\sum\limits_{Tr_k}
\frac{d^m}{q_2^m}\leqslant
\frac{2(1+2\delta)^m\lambda^m(2\varepsilon)^m}{5}.
$$
Оценим сумму по $(q_1,q_2){\in}J\setminus J_0$
$$
\sum\limits_{J\setminus J_0}\mu\left(S{\cap} A(q_1){\cap}
A(q_2)\right)\leqslant\sum\limits_{J_0}\left(\#(D{\cap}
C{\cap}\mathbb{Z}^2)\right)^m\frac{(2\varepsilon)^m}{kq_2^m}\leqslant
$$
$$
\leqslant\sum\limits_{J\setminus J_0}\left(\left(32\varepsilon
\lambda k^{1-1/m}\right)^m+ \left(\frac{16\varepsilon
k^{1-1/m}}{d}\right)^m+(8\lambda
d)^m+4^m\right)\frac{(4\varepsilon)^m}{k^{m+1}}=
$$
$$
=\frac{(128\varepsilon^2\lambda)^m}{k^2}\sum\limits_{ J\setminus
J_0}1+\frac{(64\varepsilon^2)^m}{k^2}\sum\limits_{J\setminus
J_0}\frac{1}{d^m}+\frac{(32\varepsilon\lambda)^m}{k^{m+1}}\sum\limits_{
J\setminus
J_0}d^m+\frac{(16\varepsilon)^m}{k^{m+1}}\sum\limits_{J\setminus
J_0}1\leqslant
$$
$$
\leqslant\frac{(128\varepsilon^2\lambda)^m}{k^2}\sum\limits_{
J}1+\frac{(64\varepsilon^2)^m}{k^2}\sum\limits_{J}1+\frac{(32\varepsilon
\lambda)^m}{k^{m+1}}\sum\limits_{J\setminus J_0}(9\varepsilon
k^{1-1/m})^m+\frac{(16\varepsilon
)^m}{k^{m+1}}\sum\limits_{Tr_k}1\leqslant
$$
$$
\leqslant\frac{(128\varepsilon^2
\lambda)^m}{k^2}\frac{2k^{1+\frac{3}{2m}}}{\lambda\varepsilon^2}+\frac{(64\varepsilon^2)^m
}{k^2}\frac{2k^{1+\frac{3}{2m}}}{\lambda\varepsilon^2}+\frac{(288\varepsilon^2
\lambda)^m}{k^2}\sum\limits_{J}1+\frac{(16\varepsilon
)^m}{k^{m+1}}k^2\leqslant
$$
$$
\leqslant\frac{(128\varepsilon^2
\lambda)^m}{k^2}\frac{2k^{1,75}}{\lambda\varepsilon^2}+\frac{(64\varepsilon^2)^m
}{k^2}\frac{2k^{1,75}}{\lambda\varepsilon^2}+\frac{(288\varepsilon^2
\lambda)^mk^{1,75}}{k^2}+\frac{(16\varepsilon)^m}{k}\leqslant\frac{1}{\lambda
k^{0,25}}\leqslant\frac{1}{k^{0,2}}.
$$

Найдем ограничение на сумму по $(q_1,q_2){\in}Tr_k\setminus J$
$$
\sum\limits_{Tr_k\setminus J}\mu\left(S{\cap} A(q_1){\cap}
A(q_2)\right){\leqslant}\sum\limits_{J_0}\#\left(D{\cap}
C{\cap}\mathbb{Z}^2\right)^m\frac{(2\varepsilon)^m}{kq_2^m}{\leqslant}\sum\limits_{Tr_k\setminus
J}\left(17\varepsilon\lambda
k^{1-1/m}\right)^m\frac{(2\varepsilon)^m}{k\left(\frac{k}{2}\right)^m}=
$$
$$
=\sum\limits_{Tr_k\setminus J} \frac{(68\varepsilon^2\lambda)^m
}{k^2}\leqslant\frac{(68\varepsilon^2\lambda)^m
}{k^2}\sum\limits_{Tr_k}1\leqslant\frac{(68\varepsilon^2\lambda)^m}{4}.
$$
Подставим в неравенство (\ref{ff}) неравенство (\ref{ch1}) и
неравенства, полученные выше
$$
\mu(S{\cap}\overline{M}_k)\geqslant\frac{2^m(1-\delta)^m(\lambda\varepsilon)^m}{2}-
\frac{2(1+2\delta)^m(2\lambda\varepsilon)^m}{5}-
\frac{1}{k^{0,2}}-\frac{(68\varepsilon^2\lambda)^m}{4}\geqslant
$$
$$
\geqslant\frac{2^m(1-\delta)^m(\lambda\varepsilon)^m}{2}-
\frac{2(1+3\delta)^m(2\lambda\varepsilon)^m}{5}-
\frac{(68\varepsilon^2\lambda)^m}{4}\geqslant
$$
$$
=2^m(\lambda\varepsilon)^m\left(\frac{(1-\delta)^m}{2}-
\frac{2(1+3\delta)^m}{5}\right)-\frac{(68\lambda\varepsilon^2)^m}{4}\geqslant
\frac{2^m(\lambda\varepsilon)^m}{6}-\frac{(68\lambda\varepsilon^2)^m}{4}.
$$

Лемма 14 доказана.

\textbf{4. Последовательности Бореля-Кантелли.}

Пусть здесь и далее $(\Omega,\mathcal{F},\mu)$ - вероятностное
пространство с $\Omega=[0;1]^m$.

Назовем последовательность измеримых множеств $\{M_k\}$ -
последовательностью Бореля-Кантелли если
$$
\mu\left(\bigcap\limits_{r=1}^\infty\bigcup\limits_{k=r}^\infty
M_k\right)=1.
$$
То есть почти все точки множества $\Omega$ попадают в бесконечное
количество множеств последовательности $\{M_k\}$.

{\bf Теорема Шустера \cite{Shus}.}

Пусть $\{M_k\}$ - последовательность измеримых множеств. Если для
любого измеримого множества $E$, такого что $\mu(E)>0$,
выполняется
$$
\sum\limits_{k=1}^\infty \mu(E{\cap} M_k)=+\infty,
$$
то последовательность $\{M_k\}$ является последовательностью
Бореля-Кантелли.

{\bf Теорема о плотности.} Пусть $E$ - измеримое множество с
$\mu(E)>0$, тогда для любого $\delta>0$ существует куб $S$, такой
что
$$
\frac{\mu(E{\cap} S)}{\mu(S)}>1-\delta.
$$

{\bf Лемма 15}. Пусть $\{M_k\}$ - последовательность измеримых
множеств, для любого куба $S$ со стороной $\lambda$ начиная с
некоторого номера $k_0(\lambda)$ для всех $k>k_0$ выполняется
неравенство $\mu(S{\cap} M_k)>c\mu(S)$, где $c=c(m)$, тогда
последовательность $\{M_k\}$ является последовательностью
Бореля-Кантелли.

{\bf Доказательство}. Для любых множеств $E$, $S$ и $M_k$
выполняется неравенство
\begin{equation}\label{3mn}
\mu(E{\cap} M_k)\geqslant\mu(S{\cap} E{\cap}
M_k)\geqslant\mu(S)-\mu(S{\cap}\overline{M_k})-\mu(S{\cap}\overline{E}).
\end{equation}

Из условия леммы для всех $k>k_0(\lambda)$ будет
\begin{equation}\label{u14}
\mu(S{\cap}\overline{M_k})=\mu(S)-\mu(S{\cap}
M_k)\leqslant\mu(S)-c\mu(S).
\end{equation}

По теореме о плотности для $\delta=\frac{c}{2}$ можно выбрать куб
$S$ так, что
$$
\mu(S{\cap} E)\geqslant\frac{c}{2}\mu(S),
$$
откуда получаем неравенство
\begin{equation}\label{plo}
\mu(S{\cap}\overline{E})=\mu(S)-\mu(S{\cap}
E)\leqslant\frac{c}{2}\mu(S),
\end{equation}

Подставим (\ref{u14}) и (\ref{plo}) в неравенство (\ref{3mn})
$$
\mu(E{\cap}
M_k)\geqslant\mu(S)-\mu(S)+c\mu(S)-\frac{c}{2}\mu(S)=\frac{c}{2}\mu(S).
$$

Видим, что для последовательности $\{M_k\}$ выполняются условия
теоремы Шустера, таким образом последовательность $\{M_k\}$
является последовательностью Бореля- Кантелли.

Лемма 15 доказана.

Напомним некоторые свойства прямого произведения.

Пусть $A, A'\subset X$, $B, B'\subset Y$. Тогда

1) $(A\times B){\cap}(A'\times B')=(A{\cap} A')\times(B{\cap}
B')$;

2) $\mu_{X\times Y}(A\times B)=\mu_{X}(A)\cdot\mu_{Y}(B)$.

{\bf Лемма 16.} Пусть $\{M^1_k\}$ и $\{M^2_k\}$ последовательности
в $\Omega$ и для любых кубов $S_1$ и $S_2$ начиная с некоторого
номера $k_0$ для всех $k>k_0$ выполняются неравенства
$\mu(M^1_k{\cap} S_1)>c\mu(S_1)$ и $\mu(M^2_k{\cap}
S_2)>c\mu(S_2)$, где $c=c(m)$, тогда последовательность
$\{M^1_k\times M^2_k\}$ является последовательностью
Бореля-Кантелли.

{\bf Доказательство.} Построим последовательность множеств
$\{\Psi_k\}=\left\{M^1_k\times M^2_k\right\}$.
 Покажем, что для этой
последовательности выполняются условия леммы 15. Возьмем куб
$S=S_1\times S_2$. Найдем меру пересечения
$$
\mu(S{\cap}\Psi_k)=\mu((S_1\times S_2){\cap}\left(M^1_k\times
M^2_k\right))=
$$
$$
=\mu((S_1{\cap} M^1_k)\times (S_2{\cap} M^2_k))=\mu(S_1{\cap}
M^1_k)\cdot\mu(S_2{\cap} M^2_k).
$$
По условию леммы для всех $k>k_0$ выполняется неравенство
$$
\mu(S_1{\cap} M^1_k)\cdot\mu(S_2{\cap} M^2_k)\geqslant
c\mu(S_1)\cdot c\mu(S_2)=c^2\mu(S).
$$
Видим, что для последовательности $\{\Psi_k\}$ выполняются условия
леммы 15, таким образом последовательность $\{\Psi_k\}$ является
последовательностью Бореля-Кантелли.

Лемма 16 доказана.

\textbf{5. Доказательство  основной теоремы.}

Построим последовательность множеств
$\{\Psi_k\}=\left\{\underline{M_k}\times\overline{M}_k\right\}$ и
$\{\Phi_k\}=\left\{\overline{M}_k\times\underline{M_k}\right\}$.
По леммам 3, 4, 5, 14 и 16 обе последовательности $\{\Psi_k\}$ и
$\{\Phi_k\}$ являются последовательностями Бореля-Кантелли.
Установим взаимнооднозначно соответствие между точками из $\Omega$
и матрицами.

В случае $n{=}2$ и $m{=}1$, точке $(\alpha,\beta){\in}[0,1]^2$
сопоставим матрицу
$\Theta=(\theta_1^1\,\theta_1^2)=(\alpha\,\beta)$.

В случае $n=1$ и $m\geqslant2$ точке
$(\alpha_1,\ldots,\alpha_j){\in}[0,1]^m$ сопоставим матрицу\\
$\Theta=\left(\scriptsize%
\begin{array}{c}
  \theta_1^1\\
  \vdots\\
  \theta_m^1\\
\end{array}%
\right)=
\left(\scriptsize%
\begin{array}{c}
  \alpha_1^1\\
  \vdots\\
  \alpha_m^1\\
\end{array}%
\right)$.

В обоих случаях почти все точки, а значит и почти все матрицы
попадают в последовательности  $\{\Psi_k\}$ и $\{\Phi_k\}$
бесконечное количество раз. Если $(\Theta,\Theta'){\in}\Psi_k$, то
разность $\psi_\Theta (t) - \psi_{\Theta'}(t)>0$ и если
$(\Theta,\Theta'){\in}\Phi_k$, то разность $\psi_\Theta (t) -
\psi_{\Theta'}(t)<0$, а значит для почти всех пар матриц
осцилляция происходит бесконечное количество раз.

Теорема доказана.

\end{document}